\documentstyle{amsppt}
\pagewidth{5in}
\pageheight{7.8in}
\magnification=\magstep1
\hyphenation{co-deter-min-ant co-deter-min-ants pa-ra-met-rised
pre-print pro-pa-gat-ing pro-pa-gate
fel-low-ship Cox-et-er dis-trib-ut-ive}
\def\leaderfill{\leaders\hbox to 1em{\hss.\hss}\hfill}
\def\A{{\Cal A}}
\def\D{{\Cal D}}
\def\H{{\Cal H}}

\def\idest{i.e.\ }
\def\a{{\alpha}}

\def\g{{\gamma}}

\def\d{{\delta}}
\def\De{{\Delta}}
\def\e{{\varepsilon}}

\def\k{{\kappa}}
\def\l{{\lambda}}

\def\s{{\sigma}}
\def\t{{\tau}}
\def\w{{\omega}}

\def\W{{\widehat W}}

\def\End{\text{\rm End}}
\def\Hom{\text{\rm Hom}}
\def\im{\text{\rm im}}

\def\qv{{{\Bbb Q}(v)}}
\def\ugn{{U(\widehat{{\frak g \frak l}_n})}}
\def\usn{{U(\widehat{{\frak s \frak l}_n})}}
\def\ugln{{U({\frak g \frak l}_n})}
\def\usln{{U({\frak s \frak l}_n})}
\def\sqnr{{\widehat S_q(n, r)}}
\def\boxit#1{\vbox{\hrule\hbox{\vrule \kern3pt
\vbox{\kern3pt\hbox{#1}\kern3pt}\kern3pt\vrule}\hrule}}
\def\rabbit{\vbox{\hbox{\kern0pt
\vbox{\kern0pt{\hbox{---}}\kern3.5pt}}}}

\def\tableau#1{
        \hbox {
                \hskip -10pt plus0pt minus0pt
                \raise\baselineskip\hbox{
                \offinterlineskip
                \hbox{#1}}
                \hskip0.25em
        }
}

\def\tabCol#1{
\hbox{\vtop{\hrule
\halign{\strut\vrule\hskip0.5em##\hskip0.5em\hfill\vrule\cr\lower0pt
\hbox\bgroup$#1$\egroup \cr}
\hrule
} } \hskip -10.5pt plus0pt minus0pt}

\def\CR{
        $\egroup\cr
        \noalign{\hrule}
        \lower0pt\hbox\bgroup$
}



\topmatter
\title
The affine $q$-Schur algebra
\endtitle

\author R.M. Green\endauthor
\affil 
Department of Mathematics and Statistics\\ Lancaster University\\
Lancaster LA1 4YF\\ England\\
{\it  E-mail:} r.m.green\@lancaster.ac.uk
\endaffil

\abstract
We introduce an analogue of the $q$-Schur algebra associated to Coxeter
systems of type $\widehat A_{n-1}$.  We give two constructions
of this algebra.  The first construction realizes the algebra as a certain
endomorphism algebra arising from an affine Hecke algebra of type $\widehat
A_{r-1}$, where $n \geq r$.  This generalizes the original $q$-Schur
algebra as defined by Dipper and James, and the new algebra
contains the ordinary $q$-Schur algebra and the affine Hecke algebra
as subalgebras.  Using this approach we can prove a double centralizer
property.  The second construction realizes the affine $q$-Schur
algebra as the faithful quotient of the action of a quantum 
group on the tensor power of a certain module, analogous to the
construction of the
ordinary $q$-Schur algebra as a quotient of $U(\frak g \frak l_n)$.
\endabstract


\thanks
The author was supported in part by an E.P.S.R.C. postdoctoral
research assistantship.
\endthanks
\endtopmatter

\centerline{\bf To appear in the Journal of Algebra}

\head Introduction \endhead

The $q$-Schur algebra, $S_q(n, r)$, is a finite dimensional algebra
which first appeared in the work of Dipper and James \cite{{\bf 5},
{\bf 6}} and has
applications to the representation theory of the general linear
group.  The algebra is defined in terms of $q$-permutation modules for the
Hecke algebra arising from the Weyl group of type $A$.  More recently,
analogous algebras associated to Weyl groups of type $B$ have been
studied (for example, \cite{{\bf 15}, {\bf 24}}), as well as more sophisticated
algebras of the same kind \cite{{\bf 7}, {\bf 11}} which have applications to 
representation theory.

An interesting feature of the $q$-Schur algebra (of type $A$) is that
it arises as a quotient of the quantized enveloping algebra $\ugln$.
This relationship was explored in \cite{{\bf 13}}.
It was shown in \cite{{\bf 6}} that when $n \geq r$, the $q$-Schur algebra
$S_q(n, r)$ is the centralizing partner of the action of the Hecke
algebra $\H({\Cal S}_r)$ (where $\Cal S_r$ is the symmetric group on
$r$ letters) on ``$q$-tensor space''.  Jimbo \cite{{\bf 18}} established a
quantized Weyl reciprocity between the quantized enveloping algebra
and the Hecke algebra $\H({\Cal S}_r)$ acting on a suitable tensor
space.  Combining these results means that the faithful quotient of
this action of the quantized enveloping algebra is the $q$-Schur algebra.

It is natural to wonder what the corresponding situation is for the
Weyl group of type affine $A$.  An affine analogue of Jimbo's
quantized Weyl reciprocity has been described by Chari and Pressley
\cite{{\bf 2}}.  The significant difference between the situation
considered in \cite{{\bf 2}} and the one considered in this paper is that
our tensor spaces are infinite dimensional, whereas the one in \cite{{\bf 2}}
is finite dimensional.  We wish to do this so that the tensor space
contains the regular representation of our affine Hecke algebra, which
is an infinite-dimensional algebra.  Since the representation theory
of affine Hecke algebras (see \cite{{\bf 3}, {\bf 21}}) is subtle, 
our tensor space has an interesting structure as a Hecke algebra
module.  In particular, it is not completely reducible.

The centralizer algebra of this Hecke algebra action on $q$-tensor space
is the affine $q$-Schur algebra
of the title.  This is an infinite dimensional associative algebra
which contains the affine Hecke algebra and the finite $q$-Schur
algebra as subalgebras, and like the latter algebra, it is defined in
terms of $q$-permutation modules arising from parabolic subalgebras of
a suitable Hecke algebra.  
We also exhibit two bases of this algebra: firstly a
natural one extending the Dipper--James basis of the ordinary
$q$-Schur algebra, and secondly a Kazhdan--Lusztig type basis extending Du's
basis \cite{{\bf 8}} for the $q$-Schur algebra.  

There is a rival candidate for the title of
``affine $q$-Schur algebra''; this arises (as in the case of the
ordinary $q$-Schur algebra) as a quotient of a suitable Hopf algebra
acting on a tensor power $V^{\otimes r}$ 
of a suitable natural module $V$.  We introduce a
quantum group, $\ugn$, and show that the faithful quotient of the
action gives an algebra isomorphic to the affine $q$-Schur algebra
of the preceding paragraph.  The tensor power $V^{\otimes r}$ is
isomorphic to the $q$-tensor space of the preceding paragraph,
although the isomorphism is less easy to describe than in the finite
case.  The centralizing algebra of the action $\ugn$ on
$V^{\otimes r}$ is thus isomorphic to the affine Hecke algebra, thus
establishing quantized Weyl reciprocity in the affine case.

These results raise some interesting questions about canonical bases
and other properties of $\ugn$, which we mention at the end of the
paper.

\vfill\eject
\head 1. Weyl groups and Hecke algebras \endhead

\subhead 1.1 The affine Weyl group \endsubhead

The Weyl group we consider in this paper is that of type $\widehat
A_{r-1}$, where we intend $r \geq 3$.  This corresponds to the Dynkin
diagram in Figure 1.1.

\topcaption{Figure 1.1} Dynkin diagram of type $\widehat A_{r-1}$ 
\endcaption
\centerline{
\hbox to 2.888in{
\vbox to 0.763in{\vfill
        \includegraphics{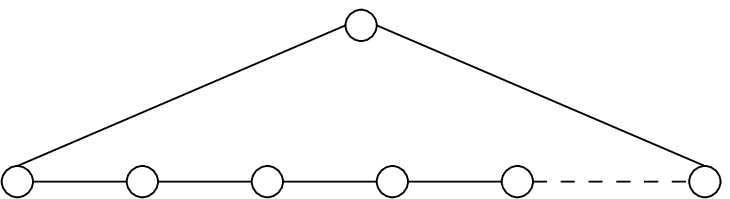}
}
\hfill}
}

The number of vertices in the graph in Figure 1.1 is $r$, as the top
vertex (numbered $r$) is regarded as an extra relative to the
remainder of the graph, which is a Coxeter graph of type $A_{r-1}$.

We associate a Weyl group, $W = W(\widehat A_{r-1})$,
to this Dynkin diagram in the usual way (as in \cite{{\bf 16}, \S2.1}).  
This associates to node $i$
of the graph a generating involution $s_i$ of $W$, where $s_i s_j =
s_j s_i$ if $i$ and $j$ are not connected in the graph, and $$
s_i s_j s_i = s_j s_i s_j
$$ if $i$ and $j$ are connected in the graph.  For $t \in {\Bbb Z}$, 
it is convenient to
denote by $\bar{t}$ the congruence class of $t$ modulo $r$, taking
values in the set $\{1, 2, \ldots, r\}$.
For the purposes of this paper, it is helpful to think of
this group as follows.

\proclaim{Proposition 1.1.1}
There exists a group isomorphism from $W$
to the set of permutations of ${\Bbb Z}$ which satisfy
the following conditions: $$
\eqalignno{
(i+r)w &= (i)w + r \text{ for } i \in {\Bbb Z}, & (1)\cr
\sum_{t = 1}^r (t)w &= \sum_{t = 1}^r t & (2)\cr
}$$
such that $s_i$ is mapped to the permutation
$$t \mapsto \cases
t & \text{ if } \bar{t} \ne \bar{i}, \overline{i + 1},\cr
t-1 & \text{ if } \bar{t} = \overline{i + 1},\cr
t+1 & \text{ if } \bar{t} = \bar{i},\cr
\endcases$$ for $t \in {\Bbb Z}$.
\endproclaim

\demo{Proof}
This is given in \cite{{\bf 22}}. 
\qed\enddemo

For reasons relating to weight spaces which will become clear later,
we consider a larger group $\W$ of permutations of ${\Bbb Z}$.

\definition{Definition 1.1.2}
Let $\rho$ be the permutation of ${\Bbb Z}$ taking $t$ to $t + 1$ for
all $t$.  Then the group $\W$ is defined to be the subgroup of
permutations of ${\Bbb Z}$ generated by the group $W$ and $\rho$.
\enddefinition

As will become clear later, the point of $\rho$ is that conjugation by
$\rho$ will correspond to a graph automorphism of the Dynkin diagram given by
rotation by one place.

\proclaim{Proposition 1.1.3}
There exists a group isomorphism from $\W$
to the set of permutations of ${\Bbb Z}$ which satisfy the
following conditions: $$
\eqalignno{
(i+r)w &= (i)w + r \text{ for } i \in {\Bbb Z}, & (1)\cr
\sum_{t = 1}^r (t)w &\equiv \sum_{t = 1}^r t \mod r.& (2)\cr
}$$
\endproclaim

\demo{Proof}
We first observe that (2) follows from (1).  (We have included (2) to
emphasise how this result is a generalization of Proposition 1.1.1.)

Next, we note that $\rho s_{i+1} \rho^{-1} = s_i$, where
$s_{r+1}$ is interpreted as meaning $s_1$.  This
shows that any element of $\W$ is expressible as $\rho^z w$ for some
integer $z$ and some $w \in W$.  Since $\rho^z$ and $w$ both have
properties (1) and (2)
given in the statement, it follows that $\rho^z w$ does also.  Thus
any element of $\W$ satisfies the properties given.

Conversely, given a permutation $w$ satisfying (1) and (2), condition
(2) shows that $$
\sum_{t = 1}^r (t)w = \left( \sum_{t = 1}^r t \right) + kr
$$ for some integer $k$.  Using the fact that $$
\sum_{t = 1}^r (t)\rho = \left( \sum_{t = 1}^r t \right) + r
,$$ we see that $w \rho^{-k}$ lies in $W$ by Proposition 1.1.1.  It
follows that $w$ lies in $\W$.
\qed\enddemo

\proclaim{Corollary 1.1.4}
Any element of $\W$ is uniquely expressible in the form $\rho^z w$ for
$z \in {\Bbb Z}$ and $w \in W$.  Conversely, any element of this form
is an element of $\W$.
\endproclaim

\demo{Proof}
The second assertion is obvious.  The first assertion follows from the
uniqueness of the integer $k$ appearing in the proof of Proposition 1.1.3.
\qed\enddemo

The group $\W$ can be expressed in a more familiar way via a
semidirect product.

\proclaim{Proposition 1.1.5}
Let $S \cong {\Cal S}_r$ be the subgroup of $\W$ generated by $$
\{s_1, s_2, \ldots, s_{r-1}\}
.$$  Let $Z$ be the
subgroup of $\W$ consisting of all permutations $z$ satisfying $$
(t)z \equiv t \mod r
$$ for all $t$.  Then ${\Bbb Z}^r \cong Z \triangleleft \W$ and $\W$
is the semidirect product of $S$ and $Z$.
\endproclaim

\demo{Proof}
It is clear that $S$ and $Z$ are isomorphic to the groups given,
since any permutation $z$ satisfying the given congruence lies 
in $\W$.  Observe
also that $S \cap Z = \{1\}$.  The other assertions follow easily.
\qed\enddemo

It is convenient to extend the usual notion of the length of an
element of a Coxeter group to the group $\W$ in the following way.

\definition{Definition 1.1.6}
If $w \in W$, the length, $\ell(w)$ of $W$ is the length of a word of
minimal length in the group generators $s_i$ of $W$ which is equal to $w$.

The length, $\ell(w')$, of a typical element $w' = \rho^z w$ of $\W$
(where $z \in {\Bbb Z}$ and $w \in W$) is defined to be $\ell(w)$.
\enddefinition

\subhead 1.2 The affine Hecke algebra \endsubhead

We now define the Hecke algebra $\H = \H(\W)$
as in \cite{{\bf 19}}.  (It will be shown in \S4.2 that this
definition is compatible with the definition given in \cite{{\bf 2}}.)
The Hecke algebra is a $q$-analogue
of the group algebra of $\W$, and is related to $W$ in
the same way as the Hecke algebra $\H({\Cal S}_r)$ of type $A$ is
related to the symmetric group ${\Cal S}_r$.  In particular, one can
recover the group algebra of $\W$ by replacing the parameter $q$
occurring in the definition of $\H(\W)$ by $1$.

\definition{Definition 1.2.1}
The affine Hecke algebra $\H = \H(\W)$ over ${\Bbb Z}[q, q^{-1}]$ is
the associative, unital algebra with algebra generators $$
\{T_{s_1}, \ldots, T_{s_r}\} \cup \{ T_\rho, T_\rho^{-1} \}
$$ and relations $$\eqalignno{
T_s^2 &= q T_s + (q-1), & (1)\cr
T_s T_t &= T_t T_s \text{ if $s$ and $t$ are not adjacent in the
Dynkin diagram}, & (2)\cr
T_s T_t T_s &= T_t T_s T_t \text{ if $s$ and $t$ are adjacent in the
Dynkin diagram}, & (3)\cr
T_\rho T_{s_{i+1}} T_\rho^{-1} &= T_{s_i}. & (4)\cr
}$$  In relation (4), we interpret $s_{r+1}$ to mean $s_1$.
\enddefinition

\definition{Definition 1.2.2}
Let $w \in W$.  The element $T_w$ of $\H(W)$ is defined as $$
T_{s_{i_1}} \cdots T_{s_{i_m}}
,$$ where $s_{i_1} \cdots s_{i_m}$ is a reduced expression for $w$
(\idest one with $m$ minimal).  (This is well-defined by standard
properties of Coxeter groups.)

If $w' \in \W$ is of form  $\rho^z w$ for $w \in W$, we denote by
$T_{w'}$ the element $$
T_\rho^z T_w
.$$  (This is well-defined by Corollary 1.1.4.)
\enddefinition

\proclaim{Proposition 1.2.3}
A free ${\Bbb Z}[q, q^{-1}]$-basis for $\H$ is given by the set $$
\{T_w : w \in \W\}
.$$
\endproclaim

\demo{Proof}
Clearly, $\H(W)$ is a quotient of the Hecke algebra arising from a
Coxeter group of type affine $A$, and so has a spanning set 
$\{T_w : w \in W\}$.

Using relation (4) of Definition 1.2.1, one can express any element of
$\H$ as a sum of elements $T_\rho^z T$, where $T$ is an element of 
the subalgebra of $\H(\W)$ generated by the elements
$T_s$.  It follows that the
set given in the statement of the proposition is a spanning set.

Suppose that there is a linear dependence
relation among the elements in the statement.  
We invoke a standard argument based on Gauss' Lemma (see  
\cite{{\bf 15}, Theorem 3.2.4}).  By clearing denominators and
negative powers of $q$, we may assume all the coefficients in the
dependence relation lie in ${\Bbb Z}[q]$ and that the largest monomial
involved is of minimal possible degree.  If the coefficients
become identically zero when $q$ is set to $1$, then $(q - 1)$ divides
each coefficient over ${\Bbb Q}(q)$ and therefore over ${\Bbb Z}[q]$.
This contradicts the minimality assumption, so the coefficients are
not identically zero when $q$ is set to $1$.  However, this gives a
dependence relation in the group algebra of $\W$, a contradiction.
\qed\enddemo

\remark{Remark 1.2.4}
The unquantized analogue of the affine Hecke algebra in the sense of
\cite{{\bf 2}, \S4.9} is precisely the semidirect product arising in 
Proposition 1.1.5.  It will turn out in \S4.2 that our algebra is the
same as the affine Hecke algebra of \cite{{\bf 2}}.
\endremark

It is sometimes useful to know that $\H$ can be generated as an
algebra more economically as follows.

\proclaim{Lemma 1.2.5}
As a ${\Bbb Z}[q, q^{-1}]$-algebra, $\H$ is generated by $T_{s_1}$, $T_\rho$
and $T_\rho^{-1}$.
\endproclaim

\demo{Proof}
This is obvious from relation (4) of Definition 1.2.1.
\qed\enddemo

\subhead 1.3 Graphs \endsubhead

It is sometimes useful to be able
to interpret the length of an element in terms of combinatoric
properties of graphs.  Any permutation of ${\Bbb Z}$ can be
represented by using the points ${\Bbb Z} \times \{0, 1\}$, where
point $i$ in the top row is connected to point $(i)w$ in the bottom
row by a straight line.

Because of condition (1) in Proposition 1.1.1 and Proposition 1.1.3,
the permutations representing elements of $\W$ have the property that
they can be represented on a cylinder with $r$ points evenly spaced on
the top and bottom faces, labelled with the congruence classes 1 up to
$r$ modulo $r$.  It is clear that one can pass between the cylindrical
and planar notation without loss.

\definition{Definition 1.3.1}
Let $w \in \W$.  We define $\nu(w)$ to be the number of crossings in the
cylindrical representation of $w$.  
\enddefinition

The advantage of representing the permutations on a cylinder can be
seen from the following result.

\proclaim{Proposition 1.3.2}
Let $w \in \W$.  Then $\ell(w) = \nu(w)$.
\endproclaim

\demo{Proof}
Note that $\nu(s_i) = 1$ and $\nu(\rho^z) = 0$.
It follows (by composition of permutations) that $\nu(w) \leq \ell(w)$
in general.

We prove the converse inequality $\nu(w) \geq \ell(w)$ by induction on
$\nu(w)$.  If $\nu(w) = 0$,
one can easily check that $(t)w = t + z$ for some $z \in {\Bbb Z}$,
and so $w = \rho^z$.  

Now consider the case where the graph of $w$ has a
crossing.  This means there exist $i < j \in {\Bbb Z}$ with $(i)w >
(j)w$.  By consideration of the pairs $(i, i+1)$, $(i+1, j-1)$ and
$(j-1, j)$, one sees that there exists such a pair $i$ and $j$ for
which $j = i + 1$.  Now the element $w s_i$ of $W'$ is such that $\nu(w
s_i) = \nu(w) - 1$.  By induction, $w s_i$ is of
length at most $\nu(w) - 1$, so $w = w s_i s_i$ is of length at most
$\nu(w)$, as desired.  This completes the induction.
\qed\enddemo

This interpretation of length leads immediately to the following corollary.

\proclaim{Corollary 1.3.3}
We have $$
(i)w < (i+1)w \Leftrightarrow \ell(s_i w) > \ell(w)
$$ and $$
(i)w^{-1} < (i+1)w^{-1} \Leftrightarrow \ell(w s_i) > \ell(w)
.$$
\endproclaim

\subhead 1.4 Distinguished coset representatives \endsubhead

A well-known property of Coxeter groups which is used extensively in
the theory of $q$-Schur algebras is that of distinguished coset
representatives with respect to parabolic subgroups.  These concepts
can be sensibly extended to the group $\W$.  Note that $\W$ is not a
Coxeter group in an obvious way, although $W$ is.

\definition{Definition 1.4.1}
Let $\Pi$ be the set of subsets of $S = \{s_1, s_2, \ldots, s_r\}$,
excluding $S$ itself.  For each $\pi \in \Pi$, we define the subgroup
$\W_\pi$ of $\W$ to be that generated by $\{s_i \in \pi\}$.  (Such a
subgroup is called a parabolic subgroup.)  We will sometimes write
$W_\pi$ for $\W_\pi$ to emphasise that it is a subgroup of $W$.

Let $\Pi'$ be the set of elements of $\Pi$ which omit the generator $s_r$.
\enddefinition

\remark{Remark 1.4.2}
All the subgroups $\W_\pi$ are in fact subgroups of $W$, and are parabolic
subgroups in the usual sense of Coxeter groups.  Furthermore, each
such $\W_\pi$ is isomorphic to a direct product of Coxeter groups of
type $A$ (\idest finite symmetric groups) corresponding to the
connected components of the Dynkin diagram obtained after omitting the
elements $s_i$ which do not occur in $\pi$.  This assertion is
immediate from well-known properties of Coxeter groups.
\endremark

\definition{Definition 1.4.3}
Let $\pi \in \Pi$.  The subset $\D_\pi$ of $\W$ is the set of
those elements such that for any $w \in \W_\pi$ and $d \in \D_\pi$, $$
\ell(w d) = \ell(w) + \ell(d)
.$$  We call $\D_\pi$ the set of distinguished right coset
representatives of $\W_\pi$ in $\W$.

The subset $\D_\pi^{-1}$ is called the set of distinguished left coset
representatives of $\W_\pi$ in $\W$;
elements $d \in \D_\pi^{-1}$ have the property that $$
\ell(d w) = \ell(d) + \ell(w)
$$ for any $w \in \W_\pi$.
\enddefinition

\proclaim{Proposition 1.4.4}
Let $\pi \in \Pi$ and $w \in \W$.  Then $w = w_\pi w^\pi$ for a unique
$w_\pi \in \W_\pi$ and $w^\pi \in \D_\pi$.  

We have $w \in \D_\pi$ if and only if 
$(t)w < (t+1)w$ whenever $(t, t+1)$ corresponds to a
generator $s_{\bar{t}}$ in $\W_\pi$.
\endproclaim

\demo{Note}
This is a familiar fact for Coxeter groups, but needs to be checked
in the case of the group $\W$.
\enddemo

\demo{Proof}
The second assertion follows from the permutation representation of
$\W$ given by Proposition 1.1.3, Corollary 1.3.3, and
the fact that $\W_\pi$ is canonically isomorphic to a direct product
of finite symmetric groups.  

Using the same techniques, we see that any $w \in \W$ arises as a
product $w_\pi w^\pi$ such that the lengths add up as claimed.
Furthermore, any element $w' w_\pi$ for $w' \in \W_\pi$ satisfies
$\ell(w' w_\pi) = \ell(w') + \ell(w_\pi)$.  It follows that the
set $\D_\pi$ is indeed a set of right coset representatives.  The
uniqueness of the decomposition also follows.
\qed\enddemo

It is clear that ${\Cal D}_\pi^{-1}$ is a set of distinguished left
coset representatives with properties analogous to the right coset
representatives given by Proposition 1.4.4.

\proclaim{Proposition 1.4.5}
Let $\pi_1, \pi_2 \in \Pi$.  The set $\D_{\pi_1, \pi_2} := 
\D_{\pi_1} \cap \D_{\pi_2}^{-1}$ is
an irredundantly described set of double
$\W_{\pi_1}$--$\W_{\pi_2}$-coset representatives, each of minimal
length in its double coset.
\endproclaim

\demo{Proof}
This follows from Proposition 1.4.4 using an elementary argument.
\qed\enddemo

The double coset representatives $\D_{\l, \mu}$ will play a key
r\^ole in the understanding the structure of the affine $q$-Schur algebra.

\head 2. The affine $q$-Schur algebra as an endomorphism algebra \endhead

\subhead 2.1 Construction of the affine $q$-Schur algebra \endsubhead

Using the affine Hecke algebra of \S1, we can now define the affine
$q$-Schur algebra.  The definition is designed to be compatible with a
quantum group which will be introduced in \S3, where the motivation
for the definitions here will become clearer.

The following definition was made in \cite{{\bf 6}}.

\definition{Definition 2.1.1}
A weight is a 
composition $\l = (\l_1, \l_2, \ldots, \l_n)$ of $r$ into $n$ pieces, 
that is, a finite sequence of nonnegative integers
whose sum is $r$.  (There is no monotonicity assumption on the
sequence.)  We denote the set of weights by $\Lambda(n, r)$.

The $r$-tuple $\ell(\l)$ of a weight $\l$ is the weakly increasing
sequence of integers where there are $\l_i$ occurrences of the entry $i$.

The Young subgroup ${\Cal S}_{\l} \subseteq {\Cal S}_r \subseteq W
\subseteq \W$ is that subgroup of permutations of the set $
\{1, 2, \ldots, r\}
$ which leaves invariant the following
sets of integers: $$
\{ 1, 2, \ldots, \l_1 \}, \{\l_1 + 1, \l_1 + 2, \ldots, \l_1 + \l_2\},
\{\l_1 + \l_2 + 1, \ldots \}, \ldots
.$$  

The weight $\w$ is given by the $n$-tuple $$
( \underbrace{1, 1, \ldots, 1}_r, \underbrace{0, 0, \ldots, 0}_{n-r} )
.$$
\enddefinition

\remark{Remark 2.1.2}
The Young subgroup ${\Cal S}_\l \subseteq {\Cal S}_r$ 
can be thought of as a group $\W_\l$ for some $\l \in \Pi'$.  Note,
however, that different compositions $\l$ can give rise to canonically
isomorphic groups.  Also note that the hypothesis $n \geq r$ is
necessary for $\w$ to exist.
\endremark

\definition{Definition 2.1.3}
Let $\l \in \Pi$.  
For $t \in {\Bbb Z}$, the parabolic subgroup $\W_{\l + t}$ is
the one generated by those elements $s_{\overline{i + t}}$ where $i$
is such that $s_i$ lies in $\W_{\l}$.
We also use the notation $\D_{\l + t}$ with the obvious meaning.

The element $x_{\l + t} \in \H$ is defined as $$
x_{\l + t} := \sum_{w \in \W_{\l + t}} T_w
.$$  We will write $x_\l$ for $x_{\l + 0}$.
\enddefinition

\definition{Definition 2.1.4}
The affine $q$-Schur algebra $\sqnr$ over ${\Bbb Z}[q, q^{-1}]$ 
is defined by $$
\sqnr := \End_\H \left( \bigoplus_{\l \in \Lambda(n, r)} x_\l \H \right),
$$ where $\H = \H(\W)$.
\enddefinition

\subhead 2.2 A basis for $\sqnr$ \endsubhead

One can describe a basis for $\sqnr$ using techniques familiar from
\cite{{\bf 4}}.  The main problem is to describe a basis for the Hom
space $\Hom_\H(x_\mu \H, x_\l \H)$.  

\proclaim{Lemma 2.2.1}
Let $\l \in \Pi'$.  As right $\H(W)$-modules, $$
x_\l \H(\W) \cong \bigoplus_{t \in {\Bbb Z}} x_\l \rho^t \H(W) \cong
\bigoplus_{t \in {\Bbb Z}} x_{\l + t} \H(W)
.$$
\endproclaim

\demo{Proof}
The first isomorphism follows by Proposition 1.2.3 and Corollary
1.1.4.  The second isomorphism follows from the definition of $x_{\l +
t}$.
\qed\enddemo

Dipper and James remark, at the
end of \cite{{\bf 4}, \S2}, that all the results of that section hold for
arbitrary finite Coxeter groups.  In fact, the results also hold for
the infinite Coxeter group $W$ with respect to finite parabolic
subgroups, for the same reasons.  This proves the following
result.

\proclaim{Lemma 2.2.2}
The set $\{\phi_{\l + t, \mu}^d : d \in W \cap \D_{\l + t, \mu}\}$ is a
${\Bbb Z}[q, q^{-1}]$-basis for $$\Hom_{\H(W)}(x_\mu \H(W), x_{\l + t}
\H(W)),$$ where $$
\phi_{\l + t, \mu}^d (x_\mu) := \sum_{d' \in {\Cal D}_\nu \cap W_\mu}
x_{\l + t} T_{d d'} = \sum_{w \in W_{\l + t} d W_\mu} T_w
$$ and $\nu$ is the composition of $n$ corresponding to the standard
Young subgroup $$d^{-1} W_{\l + t} d \cap W_\mu$$ of $W$.
\endproclaim

Using this we can deduce a basis theorem for the affine $q$-Schur
algebra.

\definition{Definition 2.2.3}
Let $d \in \W$ be an element of ${\Cal D}_{\l, \mu}$.
Write $d = \rho^z c$ (as in Corollary 1.1.4) with $c \in W$.
Then the element $$\phi_{\l, \mu}^d \in \Hom(x_\mu \H(\W), x_\l \H(\W))$$ is
defined as $$\eqalign{
\phi_{\l, \mu}^d (x_\mu) :=& 
\sum_{d' \in {\Cal D}_\nu \cap W_\mu} x_\l T_\rho^z T_{c d'}\cr
=& \sum_{d' \in {\Cal D}_\nu \cap W_\mu} T_\rho^z x_{\l + z} T_{c d'} = 
\sum_{w \in W_{\l + z} c W_\mu} T_\rho^z T_w
= \sum_{w \in W_\l d W_\mu} T_w\cr
}$$ where $\nu$ is the composition of $n$ corresponding to the standard
Young subgroup $$d^{-1} W_{\l + z} d \cap W_\mu$$ of $W$.
\enddefinition

\proclaim{Theorem 2.2.4}
A free ${\Bbb Z}[q, q^{-1}]$-basis for $\sqnr$ is given by the set $$
\{\phi_{\l, \mu}^d : \l, \mu \in \Lambda(n, r), \ d \in {\Cal D}_{\l, \mu}\}
.$$ 
\endproclaim

\demo{Proof}
The proof reduces to showing that a basis for
$\Hom_{\H(\W)}(x_\mu \H(\W), x_\l \H(\W))$ is given by $$
\{\phi_{\l, \mu}^d : d \in {\Cal D}_{\l, \mu}\}
$$  Since $x_\mu \H(\W)$ is generated by $x_\mu$, we can use lemmas
2.2.1 and 2.2.2 to reduce the problem to finding a basis for $$
\Hom_{\H(W)}(x_\mu \H(\W), x_\l \H(\W)) = 
\bigoplus_{t \in {\Bbb Z}} \Hom_{\H(W)}(x_\mu \H(W), x_{\l + t} \H(W))
.$$  Lemma 2.2.2 gives this as the set $$
\{\phi_{\l + t, \mu}^{d^\prime} : 
d \in W \cap \D_{\l + t, \mu}, t \in {\Bbb Z}\}
.$$  But since $\phi_{\l + t, \mu}^{d^\prime}$ is a homomorphism from $x_\mu$
to $x_\l T_\rho^t$, $\phi_{\l, \mu}^d$ is a homomorphism from $x_\mu$ to
$x_\l$ (set $d' = w$ and $t = z$).  Furthermore there is a natural
bijection between elements of $\D_{\l, \mu}$ and pairs $(t, W \cap
\D_{\l + t, \mu})$ for $t \in {\Bbb Z}$ given by sending $d$ to $(z,
w)$ where $d = \rho^z w$.  The proof follows.
\qed\enddemo

It is now easy to show that the affine $q$-Schur algebra contains the
affine Hecke algebra and the ordinary $q$-Schur algebra.

\proclaim{Proposition 2.2.5}
The set of basis elements $$
\{\phi_{\l, \mu}^d : \l, \mu \in \Lambda(n, r), d \in {\Cal S}_r \cap
\D_{\l, \mu}\}
$$ spans a subalgebra of $\sqnr$ canonically isomorphic to the
$q$-Schur algebra $S_q(n, r)$.

The set of basis elements $$
\{\phi_{\w, \w}^d : d \in \W\}
$$ spans a subalgebra canonically isomorphic to the Hecke algebra
$\H(\W)$, where $\phi_{\w, \w}^d$ is identified with $T_d$.
\endproclaim

\demo{Proof}
We recall that the ordinary $q$-Schur algebra is defined in formally in
the same way as $\sqnr$ (see \cite{{\bf 6}}), but with $\H$ being the
Hecke algebra of the symmetric group ${\Cal S}_r$.  The first
assertion is now immediate.

The second assertion follows from the fact that $$
\End_{\H(W)} (x_\w \H(W)) = \End_{\H(W)} (\H(W)) \cong \H(W)
,$$ since $\H(W)$ is an associative algebra with $1$.
\qed\enddemo

\subhead 2.3 $q$-tensor space and the double centralizer property \endsubhead

The definition of the $q$-Schur algebra in \S2.2 allows the notion of
$q$-tensor space to be introduced in exactly the same way as for the
finite $q$-Schur algebra in \cite{{\bf 6}}.  The motivation for the term
``tensor'' will become clearer in \S3.

\definition{Definition 2.3.1}
Identifying $T_d \in \H(\W)$ with $\phi_{\w, \w}^d \in \sqnr$, we
identify the basis element $\phi_{\l, w}^d$ of $\sqnr$ with the
coset $x_\l T_d$.  This makes $$
\bigoplus_{\l \in \Lambda(n, r)} x_\l \H = \langle \phi_{\l, w}^d : \l
\in \Lambda(n, r), d \in {\Cal D}_\l \rangle
$$ into a $\sqnr$--$\H(\W)$ bimodule, which we call $q$-tensor space,
$E(n, r)$.
\enddefinition

\proclaim{Lemma 2.3.2}
As a left $\sqnr$-module, $q$-tensor space is generated by $\phi_{\w,
\w}^1$.
\endproclaim

\demo{Proof}
This is immediate from the fact that $\phi_{\l, \w}^d . \phi_{\w,
\w}^1 = \phi_{\l, \w}^d$.
\qed\enddemo

Using the notion of $q$-tensor space, we can prove that the actions of
the affine $q$-Schur algebra and the affine Hecke algebra on
$q$-tensor space are in reciprocity with each other, that is, that
each is the centralizing algebra of the other.  This is a
straightforward generalization of \cite{{\bf 6}, Theorem 6.6}.

\proclaim{Theorem 2.3.3 (the double centralizer property)}
\item{\rm (i)}
{$\End_{\sqnr}(E(n, r)) = \H(\W).$}
\item{\rm (ii)}
{$\End_{\H(\W)}(E(n, r)) = \sqnr.$}
\endproclaim

\demo{Proof}
The proof of (ii) follows easily from Definition 2.3.1 and the
definition of $\sqnr$.  

To prove (i), first observe that  
Lemma 2.3.2 implies that any endomorphism in $\End_{\sqnr}(E(n, r))$ is
determined by its effect on the element $x_\w$ in the subspace
$x_\w \H$ of $q$-tensor space.  Since $$
\phi_{\l, \l}^1 \phi_{\mu, \w}^d = \delta_{\l, \mu} \phi_{\mu, w}^d
,$$ this endomorphism must map $x_\w \H$ to itself.

The proof of Proposition 2.2.5 shows that these endomorphisms are in
canonical correspondence with elements of $\H(W)$ acting by right
multiplication.
\qed\enddemo

\subhead 2.4 The Kazhdan--Lusztig basis \endsubhead

Using the methods of Du \cite{{\bf 8}}, we can describe a second basis for
the algebra $\sqnr$, analogous to the Kazhdan--Lusztig basis for the
Hecke algebra of a Coxeter group \cite{{\bf 20}}.  Such a basis is called
an IC (intersection cohomology) basis; this terminology is justified
in the case of the ordinary $q$-Schur algebra, where the basis has intimate
connections to Lusztig's canonical basis for $U^+$ (the ``plus part''
of a quantized enveloping algebra) and the theory of perverse sheaves (see
\cite{{\bf 14}, Corollary 4.7}).

The following result is well known for Coxeter groups and their
finite parabolic subgroups (each of which has a unique
longest element).

\proclaim{Lemma 2.4.1}
Let $\l, \mu \in \Pi$ and $w \in \W$.  The double coset $\W_\l w \W_\mu$
has a unique element of maximal length, denoted by $w^+$.
\endproclaim

\demo{Proof}
We may assume that $w \in {\Cal D}_{\l, \mu}$.  Write $w = \rho^t w'$
for $w' \in W$.  Then $$
\W_\l w \W_\mu
= \rho^t \W_{\l + t} w' \W_\mu
.$$  Using the result for Coxeter groups, we know that $\W_{\l + t} w'
\W_\mu$ has a longest element, $w^{\prime +}$.  It is easy to check
that the element $w^+ = \rho^t w^{\prime +}$ satisfies the required criteria.
\qed\enddemo

\definition{Definition 2.4.2}
We denote by $\leq$ the (strong) Bruhat order on a Coxeter group (see
\cite{{\bf 16}, \S5.9} for the definition).  
We extend this order to the group $\W$ by stipulating that $\rho^a
y \leq \rho^b w$ (for $y, w \in W$) if and only if $a = b$ and $y \leq
w.$

Let $P_{y, w}$ be the Kazhdan--Lusztig polynomial for $y, w \in
W$.  This polynomial is guaranteed to be zero unless $y \leq w$.  
We extend these polynomials to the group $W$ by defining $P_{\rho^a y,
\rho^b w} := \d_{a, b} P_{y, w}$.

The indeterminate $v$ is defined to be $q^{1/2}$.

The element $w_{0, \mu}$ is the longest element of $\W_\mu$.
\enddefinition

We can now define the Kazhdan--Lusztig basis of $\sqnr$, taken over
the base ring $\A := {\Bbb Z}[v, v^{-1}]$.

\definition{Definition 2.4.3}
Let $\{ \phi_{\l, \mu}^d \}$ be the basis of $\sqnr$ defined in
\S2.2.  We define $$
\theta_{\l, \mu}^d : = v^{\ell(w_{0, \mu})} \sum_{z \in {\Cal D}_{\l,
\mu}} \a_{z, w} \phi_{\l, \mu}^d,
$$ where $$\a_{z, w} : = v^{-\ell(w^+)} P_{z^+, w^+}.$$
\enddefinition

\remark{Remark 2.4.4}
Note that the sum occurring above is a finite sum, because for any $w
\in \W$, there is only a finite number of elements $y$ such that $y
\leq w$.

It is clear that the basis $\{\theta_{\l, \mu}^d\}$ is a free
$\A$-basis for $\sqnr$ since the basis change matrix from the basis
$\{\phi_{\l, \mu}^d\}$ is a unit multiple of a unitriangular matrix
with respect to a suitable ordering.
\endremark

\subhead 2.5 Generators and relations \endsubhead

For the purposes of $\S3$, it is convenient to have a presentation of
$\sqnr$ over $\qv$ via generators and relations.

\proclaim{Proposition 2.5.1}
The algebra $\qv \otimes \sqnr$ is generated by elements $$
\{ \phi_{\w, \w}^d : d \in \W \} 
\cup \{ \phi_{\l, \w}^1 : \l \in \Lambda(n, r) \}
\cup \{ \phi_{\w, \l}^1 : \l \in \Lambda(n, r) \}
.$$  The elements $\phi_{\w, \w}^d$ are subject to the relations of
the affine Hecke algebra of Definition 1.2.1 under the identification
given by Proposition 2.2.5.  The generators are also subject to the
following relations, where $s$ denotes a generator 
$s_i \in \W_\l$. $$\eqalignno{
\phi_{\w, \l}^1 \phi_{\mu, \w}^1 &= \d_{\l, \mu} \sum_{d \in \W_\l}
\phi_{\w, \w}^d, & (1)\cr
\phi_{\w, \w}^s \phi_{\w, \l}^1 &= q \phi_{\w, \l}^1, & (2) \cr
\phi_{\l, \w}^1 \phi_{\w, \w}^s &= q \phi_{\l, \w}^1. & (3) \cr
}$$
\endproclaim

\demo{Note}
Note that $s \ne s_r$ above, since $\l \in \Pi'$.
\qed\enddemo

\demo{Proof}
Using Proposition 2.2.5 and properties of the ordinary $q$-Schur
algebra (or the definition of the $\phi$-basis), it can be verified
that the relations given are all true in $\sqnr$.


Let $d = \rho^z c \in \D_{\l, \mu}$, where $c \in W$ (and thus $c \in
\D_{\l + z, \mu}$).  It follows from Definition 2.2.3 that $$\eqalign{
\phi_{\l, \w}^1 \phi_{\w, \w}^d \phi_{\w, \mu}^1 (x_\mu) &= 
T_\rho^z x_{\l + z} T_c x_\mu\cr
&= \left( \sum_{d' \in {\Cal D}_\nu \cap W_\mu} T_\rho^z x_{\l + z} T_{c
d'} \right) x_\nu\cr
&= P_\nu(q) \sum_{d' \in {\Cal D}_\nu \cap W_\mu} T_\rho^z x_{\l + z} T_{c
d'}\cr
&= P_\nu(q) . \phi_{\l, \mu}^d
,}$$ where $\nu$ is as in Definition 2.2.3 and $P_\nu$ is the Poincar\'e
polynomial $$
P_\nu(q) := \sum_{w \in \W_{\nu}} q^{\ell(w)}
.$$  This relies on the fact that the expression in parentheses in the
second line corresponds to a union of left cosets of $\W_\mu$, and
$\W_\nu \subseteq \W_\mu$.

It follows that the given set does generate $\sqnr$, since
$P_\nu(q)$ is invertible.  It also follows that the relations given
are sufficient to express the product of two basis elements as a
linear combination of others, which completes the proof.
\qed\enddemo

\head 3. The affine $q$-Schur algebra as a Hopf algebra quotient \endhead

The aim of \S3 is to realize the affine $q$-Schur algebra as the
faithful quotient of the action of a certain quantum group acting on a
tensor power of a natural module.  This explains why the definition of
the $q$-Schur algebra as given in \S2 is a natural one.  The tensor
power of the natural module will be seen to be the analogue of the
$q$-tensor space of \S2, thus justifying the terminology.

In this section, we take the base ring of $\sqnr$ and $\H(\W)$
to be $\qv$ by tensoring in the obvious way.

\subhead 3.1 The Hopf algebra $\ugn$ \endsubhead

We now introduce a quantum group, $\ugn$, based on the quantized
enveloping algebra $\usn$ associated to a Dynkin diagram of type
$\widehat A_{n-1}$.  From now on, $\bar{\ }$ will denote congruence
modulo $n$ (as opposed to modulo $r$), unless otherwise specified.

The main differences between $\ugn$ and $\usn$ are that the element $K_i$
in $\usn$ corresponds to the element $K_i K_{\overline{i+1}}^{-1}$ in
$\ugn$ (similar to the difference between
$\ugln$ and $\usln$), and, more strikingly, $\ugn$ contains a new
grouplike element, $R$, which will play an important r\^ole in the
results to come.

\definition{Definition 3.1.1}
The associative, unital algebra $\ugn$ over $\qv$ is given by
generators $$E_i, F_i, K_i, K_i^{-1}, R, R^{-1}$$ 
(where $1 \leq i \leq n$) subject to the following relations: 
\vfill\eject 
$$\eqalignno{
&K_i K_j = K_j K_i, & (1)\cr
&K_i K_i^{-1} = K_i^{-1} K_i = 1, & (2)\cr
&K_i E_j = v^{{\epsilon^{+}}(i, j)} E_j K_i, & (3)\cr
&K_i F_j = v^{{\epsilon^{-}}(i, j)} F_j K_i, & (4)\cr
&E_i F_j - F_j E_i = \delta_{ij} {
{K_i K_{i + 1}^{-1} - K_i ^{-1} K_{i+1}} \over {v - v^{-1}}}, & (5)\cr
&E_i E_j = E_j E_i \quad \text{ if $i$ and $j$ are not adjacent,} & (6)\cr
&F_i F_j = F_j F_i \quad \text{ if $i$ and $j$ are not adjacent,} & (7)\cr
&E_i^2 E_j - (v + v^{-1}) E_i E_j E_i + E_j E_i^2 = 0 \quad 
\text{ if $i$ and $j$ are adjacent,} & (8)\cr
&F_j^2 F_i - (v + v^{-1}) F_j F_i F_j + F_i F_j^2 = 0 \quad
\text{ if $i$ and $j$ are adjacent,} & (9)\cr
&R R^{-1} = R^{-1} R = 1, & (10)\cr
&R^{-1} K_{\overline{i+1}} R = K_i, & (11)\cr
&R^{-1} K_{\overline{i+1}}^{-1} R = K_i^{-1}, & (12)\cr
&R^{-1} E_{\overline{i+1}} R = E_i, & (13)\cr
&R^{-1} F_{\overline{i+1}} R = F_i. & (14)\cr
}$$ By ``adjacent'', we mean that $i$ and $j$ index adjacent nodes in
the Dynkin diagram.  Here,
$$\epsilon^+ (i, j) := \cases
1 & \text{ if } j = i;\cr
-1 & \text{ if } \bar{j} = \overline{i-1};\cr
0 & \text{ otherwise;}\cr
\endcases $$ and
$$\epsilon^- (i, j) := \cases
1 & \text{ if } \bar{j} = \overline{i-1};\cr
-1 & \text{ if } j = i;\cr
0 & \text{ otherwise.}\cr
\endcases$$
\enddefinition

\remark{Remark 3.1.2}
Although it appears at first that the set of defining relations is
rather large, it should be noted that in fact the algebra may be
generated using only $E_1, F_1, K_1, K_1^{-1}, R, R^{-1}$, since the
other generators may be expressed in terms of these by using relations
(11)--(14).  (Compare this with Lemma 1.2.5.)  Also note that $R^n$
lies in the centre of $\ugn$, just as $T_\rho^r$ lies in the centre of
$\H(\W)$.
\endremark

It will turn out that the algebra $\ugn$ can be given the structure of
a Hopf algebra (which justifies the use of the term ``quantum
group'').  This is an extension of the well-known Hopf algebra
structure \cite{{\bf 17}, \S4.8} of the quantized enveloping algebras associated to
Kac--Moody Lie algebras.  Since many of the arguments are familiar, we
will highlight only the differences between the situation of $U = \ugn$
and quantized enveloping algebras.
Although we do not need the antipode in the Hopf algebra for our
purposes, we introduce it for completeness.

\proclaim{Lemma 3.1.3}
There is a unique homomorphism of unital $\qv$-algebras $\De : U
\rightarrow U \otimes U$ such that $$\eqalign{
\De(1) &= 1 \otimes 1,\cr
\De(E_i) &= E_i \otimes K_i K_{i+1}^{-1} + 1 \otimes E_i,\cr
\De(F_i) &= K_i^{-1} K_{i+1} \otimes F_i + F_i \otimes 1,\cr
\De(X) &= X \otimes X \text{\rm \ for } X \in \{K_i, K_i^{-1}, R, R^{-1}\}.\cr
}$$
\endproclaim

\demo{Note}
The elements $X$ above are known as grouplike elements.
\enddemo

\demo{Proof}
This follows by a routine check, which is familiar except for the
cases involving the elements $R$ and $R^{-1}$.
\qed\enddemo

\proclaim{Lemma 3.1.4}
The map $\De$ is coassociative.  In other words, $(\De \otimes 1)\De =
(1 \otimes \De)\De$.
\endproclaim

\demo{Proof}
Since $\De$ is an algebra homomorphism, it suffices to check this on
generators.  The only unfamiliar case involves $R$ and $R^{-1}$; in
this case, both sides send $R$ to $R \otimes R \otimes R$ and $R^{-1}$
to $R^{-1} \otimes R^{-1} \otimes R^{-1}$.
\qed\enddemo

\proclaim{Lemma 3.1.5}
There is a unique homomorphism $\e : U \rightarrow \qv$ 
of unital $\qv$-algebras such that $$
\e(E_i) = \e(F_i) = 0
$$ and $$
\e(K_i) = \e(K_i^{-1}) = \e(R) = \e(R^{-1}) = 1
.$$
\endproclaim

\demo{Proof}
It is trivial to check that $\e$ preserves all the relations.
\qed\enddemo

\proclaim{Lemma 3.1.6}
For any $u \in U$, we have $$
(\e \otimes 1)\De(u) = 1 \otimes u
$$ and $$
(1 \otimes \e)\De(u) = u \otimes 1
.$$
\endproclaim

\demo{Proof}
It is enough to verify these claims on the generators.  This is easy
(including the new cases of $u \in \{R, R^{-1}\}$).
\qed\enddemo

\proclaim{Lemma 3.1.7}
There is a $\qv$-linear map $S : U \rightarrow U$ uniquely determined
by the following properties: $$\eqalign{
S(E_i) &= -E_i K_i^{-1} K_{i+1},\cr
S(F_i) &= -K_i K_{i+1}^{-1} F_i,\cr
S(K_i) &= K_i^{-1},\cr
S(K_i^{-1}) &= K_i,\cr
S(R) &= R^{-1},\cr
S(R^{-1}) &= R,\cr
\forall a, b \in U, \ S(ab) &= S(b) S(a).\cr
}$$
\endproclaim

\demo{Proof}
We verify that $S$ preserves the relations of $U$, from which it
follows that it is determined by its effect on the generators.  The
only new relations to be checked are relation (10) of Definition 3.1.1
(which is clear) and relations (11)--(14), which are of the form $$
S(R^{-1} X_{i+1} R) = S(X_i)
$$ for $X_i \in \{E_i, F_i, K_i, K_i^{-1}\}$.  Expanding the left hand
side, we obtain $$
S(R^{-1} X_{i+1} R) = S(R) S(X_{i+1}) S(R^{-1}) = R^{-1} S(X_{i+1}) R
.$$  It is then easily verified that $R^{-1} S(X_{i+1}) R = S(X_i)$ in
each case.
\qed\enddemo

\definition{Definition 3.1.8}
The $\qv$-algebra homomorphism $\eta : \qv \rightarrow U$ is defined
to send $1 \in \qv$ to $1 \in U$.

We denote the multiplication map in $U$ by $\mu$.
\enddefinition

\proclaim{Lemma 3.1.9}
We have $$
\mu(S \otimes 1)\De = \eta \e
$$ and $$
\mu(1 \otimes S)\De = \eta \e
.$$
\endproclaim

\demo{Proof}
We first check these assertions on the generators.  The only new cases
involve the elements $R$ and $R^{-1}$, which involve straightforward
checks such as $$
\mu(S \otimes 1) \De(R) = \mu(S \otimes 1)(R \otimes R) = \mu(R^{-1}
\otimes R) = 1 = \eta(1) = \eta \e(R)
.$$  This does not complete the proof, because the map $S$ is not an
algebra homomorphism, but there is a standard argument (given in
\cite{{\bf 17}, \S3.7}) which shows
that if the relations in the statement hold on the generators, they
hold in general.
\qed\enddemo

\proclaim{Theorem 3.1.10}
The algebra $\ugn$ is a Hopf algebra with multiplication $\mu$, unit
$\eta$, comultiplication $\De$, counit $\e$ and antipode $S$.
\endproclaim

\demo{Proof}
This is immediate from lemmas 3.1.3, 3.1.4, 3.1.5, 3.1.6, 3.1.7 and 3.1.9.
\qed\enddemo

\subhead 3.2 Tensor space \endsubhead

Let $V$ be the $\qv$-vector space with basis $\{e_t : t \in {\Bbb Z}\}$.
This has a natural $\ugn$-module structure as follows. 

\proclaim{Lemma 3.2.1}
There is a left action of $\ugn$ on $V$ defined by the conditions $$\eqalign{
E_i e_{t+1} &= e_t 
\text{ if } i = t \mod n, \cr
E_i e_{t+1} &= 0 
\text{ if } i \ne t \mod n, \cr
F_i e_{t} &= e_{t+1}
\text{ if } i = t \mod n, \cr
F_i e_{t} &= 0 
\text{ if } i \ne t \mod n, \cr
K_i e_t &= v e_t
\text{ if } i = t \mod n, \cr
K_i e_t &= e_t
\text{ if } i \ne t \mod n, \cr
R e_t &= e_{t+1}.\cr
}$$
\endproclaim

\demo{Proof}
The actions of $K_i^{-1}$ and $R^{-1}$ are clear from the above
conditions and relations (2) and (10) of Definition 3.1.1.
One then checks that the other relations are preserved.
\qed\enddemo

Since $\ugn$ is a Hopf algebra, the tensor product of two
$\ugn$-modules has a natural $\ugn$-module structure via the
comultiplication $\De$.

\definition{Definition 3.2.2}
The vector space $V^{\otimes r}$ has a natural $\ugn$-module structure
given by $u . v = \De(u)^{(r-1)} . v$.  We call this module tensor space.
\enddefinition

\remark{Remark 3.2.3}
The definition of $\sqnr$ in \S2 is chosen so that the resulting
$q$-tensor space is compatible with the tensor space above.  This
relationship will be explored later.
\endremark

It is convenient to introduce certain elements of $\End_U(V^{\otimes
r})$, which we call $y_i$ since they are
reminiscent of elements of the Hecke algebra in the sense of \cite{{\bf 2},
\S3.1}.  This similarity will be explored in \S4.2.

\definition{Definition 3.2.4}
For each $1 \leq t \leq r$, we define the invertible linear map $y_t$ on
$V^{\otimes r}$ by the condition $$
(e_{j_1} \otimes \cdots \otimes e_{j_r}).y_t := 
e_{j_1} \otimes \cdots \otimes e_{j_{t-1}} \otimes e_{-n + j_t} \otimes
e_{j_{t+1}} \otimes \cdots \otimes e_{j_r}.
$$
\enddefinition

We denote the subspace of $V$ spanned by the elements $\{e_1, \ldots,
e_n\}$ by $V_n$.  (This is the vector space which plays the r\^ole of
$V$ in the case of the ordinary $q$-Schur algebra.)  The following
result is a useful tool for reducing problems about $\sqnr$ to that of
the finite $q$-Schur algebra.

\proclaim{Proposition 3.2.5}
The elements $y_i$, $y_i^{-1}$ lie in $\End_U(V^{\otimes r})$.  Therefore the
action of an element $u \in \ugn$ on tensor space is determined by its
action on $V_n^{\otimes r}$.
\endproclaim

\demo{Proof}
Since $R^n$ lies in the centre of $\ugn$, the linear maps which send
$e_t$ to $e_{t \pm n}$ lie in $\End_U(V)$.  Using the comultiplication,
this result can be extended to $\End_U(V^{\otimes r})$, proving the
first assertion.

The second assertion holds because $V^{\otimes r} = V_n^{\otimes r}
. Y$, where $Y$ is the algebra of endomorphisms of $V^{\otimes r}$
generated by the $y_i$ and their inverses.
\qed\enddemo

We require the concept of a weight space of tensor space for the
results which follow.  This is a straightforward generalization of the
situation in finite type $A$.

\definition{Definition 3.2.6}
The weight $\l \in \Lambda(n, r)$ of a basis element $$
e_{t_1} \otimes e_{t_2} \otimes \cdots \otimes e_{t_r}
$$ of $V^{\otimes r}$ is given by the condition $$
\l_i := | \{ j : t_j \equiv i \mod n \}|.$$

The $\l$-weight space, $V_\l$, of $V^{\otimes r}$ is the span of all the basis
vectors of weight $\l$.
\enddefinition

\subhead 3.3 Relationship between $\ugn$ and $\H(\W)$ \endsubhead

As one might expect from the situation in \S2, the weight space $V_\w$
is of particular importance.  In this section we define certain operations
on this space arising from the action of elements of $\ugn$.  The
motivation for doing this is that these operations will generate an
algebra isomorphic to $\H(\W)$.  This will provide a link with the
results of \S2.

\definition{Definition 3.3.1}
For each $1 \leq i \leq r$, let 
$\t(T_{s_i}) : V_\w \rightarrow V_\w$
be the endomorphism corresponding to the
action of $v F_i E_i - 1 \in \ugn$.  Similarly, let $\t(T_{\rho^{-1}})
$ be the endomorphism corresponding to $$
F_n F_{n-1} \cdots F_{r+1} R
,$$ and let $\t(T_\rho)$ be the endomorphism corresponding
to $$
E_r E_{r+1} \cdots E_{n-1} R^{-1}
.$$
\enddefinition

\proclaim{Lemma 3.3.2}
The endomorphisms $\t(T_w)$ defined above satisfy the relations of
Definition 1.2.1 (after replacing $T_w$ by $\t(T_w)$).
\endproclaim

\demo{Proof}
Using the epimorphism $\a_r : \ugln \twoheadrightarrow S_q(n, r)$
described in \cite{{\bf 13}}, one finds that the action of $\t(T_{s_i})$ on
$V_\w$ in the case where $i \ne r$ corresponds to the action of
$\phi_{\w, \w}^{s_i} \in S_q(n, r)$.  (Recall from \cite{{\bf 13}} 
that $S_q(n, r)$ is the quotient of $\ugln$ by the annihilator of 
$V_n^{\otimes r}$.)  This proves relations (1)--(3) of Definition
1.2.1 in the case where $s_r$ is not involved.

It suffices now to prove relation (4), from which the other cases of
relations (1)--(3) may be deduced.  
The effect of $\t(T_\rho)$ on
$V_\w$ is $$
\t(T_\rho)(e_{i_1} \otimes \cdots \otimes e_{i_r})
= e_{j_1} \otimes \cdots \otimes e_{j_1}
,$$ where $j_t = i_t - 1 \mod r$.
The effect of $\t(T_{\rho^{-1}})$ on
$V_\w$ is the inverse of this action.  The proof of relation (4) now follows by
calculation of the action of $v F_i E_i - 1$ on $V_\w$ using the
comultiplication.
\qed\enddemo

In the following result, we temporarily
replace the $v$ occurring in the definition of the action of
$\t(T_w)$ by $1$, in order to regard the endomorphisms $\t(T_w)$ as
endomorphisms of a version of tensor space over ${\Bbb Q}$.
This modified tensor space is the $r$-fold tensor power of the 
${\Bbb Q}$-vector space with basis $\{e_t : t \in {\Bbb Z}\}$.  We
call this process specializing $v$ to $1$.

\proclaim{Lemma 3.3.3}
Specialize the parameter $v$ to $1$.
Then the map taking $\t(T_w)$ to $w$ extends uniquely to an 
isomorphism of algebras between the group algebra ${\Bbb Q}\W$ of
$\W$ and
the algebra $\t(\H)$ generated by the
endomorphisms $\t(T_w)$ of Definition 3.3.1.
\endproclaim

\demo{Proof}
It follows from Lemma 3.3.2 that the algebra generated by the
$\t(T_w)$ is a quotient of the group algebra of $\W$.  

To complete the proof, it is sufficient (by Proposition 1.1.3) to show
that $w \in \W$ corresponds under the isomorphism given 
to an endomorphism $\t \in \t(\H)$ such that $$
\t(e_{i_1} \otimes \cdots \otimes e_{i_r})
=
e_{w(i_1)} \otimes \cdots \otimes e_{w(i_r)}
.$$ (The left action of $\W$ is defined by
$\rho(t) = t - 1$ and $s_1$ acting via the simple transposition $(1,
2)$.  It follows from Proposition 1.1.3 that this is faithful.)  

This claim is easily checked when $w$ is one of the generators 
$s_1$, $\rho$ or $\rho^{-1}$ of $\W$.  The general case follows by
analysing the action of $\t(T_w)$ on $V_\w$ in the case $v = 1$.  A
check shows that this is given by $$
\t(T_w)(e_{i_1} \otimes \cdots \otimes e_{i_r})
=
e_{w(i_1)} \otimes \cdots \otimes e_{w(i_r)}
.$$  The proof follows from this.
\qed\enddemo

We can obtain an analogue of this result which holds in the quantized case.

\proclaim{Lemma 3.3.4}
The map taking $\t(T_w)$ to $T_w$ extends uniquely to an 
isomorphism of algebras between $\H(\W)$ and 
the algebra $\t(\H)$ generated by the
endomorphisms $\t(T_w)$ of Definition 3.3.1.
\endproclaim

\demo{Proof}
Lemma 3.3.2 shows that the elements $\t(T_w)$ satisfy the required
defining relations, so that $\t(T_w)$ is a quotient of $\H(\W)$.

Suppose an element of $\H(\W)$ maps to zero under this quotient map.
Then by the usual technique of clearing denominators and dividing out
by $v - 1$, we can find an element of ${\Bbb Q}\W$ which maps to zero
under the isomorphism of Lemma 3.3.3.  This completes the proof.
\qed\enddemo

\subhead 3.4 The map $\a_r : \ugn \twoheadrightarrow \sqnr$ \endsubhead

The results of \S3.3 can be extended to the whole tensor space
$V^{\otimes r}$ by using the generators and relations for $\sqnr$
which were given in \S2.5.  This construction relies on the surjective
homomorphism $\a_r$ from $\ugln$ to $S_q(n, r)$ which was studied in 
\cite{{\bf 13}}.

\definition{Definition 3.4.1}
For each $\l \in \Lambda(n, r)$, choose $$\eqalign{
E(\l) & \in \a_r^{-1}(\phi_{\l, \w}^1),\cr
F(\l) & \in \a_r^{-1}(\phi_{\w, \l}^1),\cr
G(\l) & \in \a_r^{-1}(\phi_{\l, \l}^1).\cr
}$$  We also denote by $E(\l)$
(respectively, $F(\l)$, $G(\l)$) the elements of $\ugn$ which are the images of
$E(\l)$ (respectively, $F(\l)$, $G(\l)$)
under the obvious homomorphism from $\ugln$ to $\ugn$.

Denote by $\g$ the map from $\ugn$ to $\End(V^{\otimes r})$
corresponding to the action on tensor space arising from the comultiplication.
\enddefinition

\proclaim{Lemma 3.4.2}
There is a unique homomorphism of algebras $$\k : \sqnr \rightarrow 
\End(V^{\otimes r})$$ such that for all $\l \in \Lambda(n, r)$, $$\eqalign{
\k(\phi_{\l, \w}^1) &= \g(E(\l)),\cr
\k(\phi_{\w, \l}^1) &= \g(F(\l)),\cr
\k(\phi_{\l, \l}^1) &= \g(G(\l)),\cr
\k(\phi_{\w, \w}^{s_1}) &= \t(T_{s_1})\g(G(\w)),\cr
\k(\phi_{\w, \w}^{\rho^{\pm 1}}) &= \t(T_{\rho^{\pm 1}})\g(G(\w)).\cr
}$$
\endproclaim

\demo{Proof}
If such a homomorphism exists, it is unique by Lemma 1.2.5 and
Proposition 2.5.1.  

It remains to check the relations of Proposition 2.5.1.  The Hecke
algebra style relations hold because of Lemma 3.3.4.  The other
relations are essentially relations of the ordinary $q$-Schur
algebra.  These hold because of Proposition 3.2.5.
\qed\enddemo

\proclaim{Lemma 3.4.3}
Let $w \in {\Cal S}_r$, $\l \in \Lambda(n, r)$, and let $w = w_\l d$ be the
decomposition of $w$ as in Proposition 1.4.4, so that $d \in \D_\l$.
Consider the basis element $$
e_{i_1} \otimes \cdots \otimes e_{i_r} = 
e_{(1)w^{-1}} \otimes \cdots \otimes e_{(r)w^{-1}}
$$ of $V_n^{\otimes r} \subset V^{\otimes r}$.  Then the action of
$E(\l)$ sends this element to one of the form $$
v^{2 \ell(w_\l) + f(\l)} e_{j_1} \otimes \cdots \otimes e_{j_r}
,$$ where $j_t$ is the $i_t$-th component of $\ell(\l)$ as in
Definition 2.1.1, and $f$ is a certain integer-valued function.
\endproclaim

\demo{Proof}
This is essentially a statement about the ordinary $q$-Schur algebra,
and is therefore a consequence of \cite{{\bf 12}, Theorem 3.8}.  Alternatively,
one can argue as follows.

The claim holds in the case $w = 1$ using the definition of $q$-tensor space
(see \cite{{\bf 6}, Definition 2.6} and the discussion following it), and
the fact that $q$-tensor space and tensor space are canonically
isomorphic (in the finite dimensional case) by multiplication by a
suitable integer power of $v$ (see \cite{{\bf 9}, \S1.5}).  Note that this
canonical isomorphism respects the notion of weight space
(see Definition 2.1.1 and Definition 3.2.6).

The general case uses the fact that the actions of $\H({\Cal S}_r)$
and $\ugln$ on $V_n^{\otimes r}$ commute, and the definition of the
action of $\H({\Cal S}_r)$ as given in \cite{{\bf 9}}.
\qed\enddemo

\remark{Remark 3.4.4}
This result determines the action of $E(\l)$ on $V_\w$ because of
Proposition 3.2.5.
\endremark

\proclaim{Lemma 3.4.5}
Consider the basis element $$
e_{\ell(\mu)_1} \otimes \cdots \otimes e_{\ell(\mu)_r}
$$ of $V_n^{\otimes r} \subset V^{\otimes r}$.  Then the action of
$F(\mu)$ sends this element to one of the form $$
v^{g(\mu)} \sum_{w \in \W_\mu} v^{\ell(w)}
e_{(1)w^{-1}} \otimes \cdots \otimes e_{(r)w^{-1}}
,$$ where $g$ is a certain integer-valued function.
\endproclaim

\demo{Proof}
This is proved using the same techniques as in the proof of Lemma
3.4.3, although this is an easier case.
\qed\enddemo

\proclaim{Lemma 3.4.6}
We have $\g(U) = \im(\k)$.  Thus $\im(\k)$ is a quotient of
$U$ by the kernel of its action on tensor space.
\endproclaim

\demo{Proof}
The second assertion follows if we can establish $\g(U) \subseteq
\im(\k)$, because $\k(\phi)$ gives the same endomorphism as the
action of a suitable element of $\ugn$ for $\phi$ a generator of $\sqnr$.

For the first assertion, it 
is enough to prove that $\g(u) \in \im(\k)$ for $$u \in \{E_1, F_1,
K_1, K_1^{-1}, R, R^{-1}\}.$$  The case where $u \in
\{E_1, F_1, K_1, K_1^{-1}\}$ holds because of the existence of the
surjection $\a_r$ in the finite case and Proposition 3.2.5.

We now show that $\g(R) \in \im(\k)$.  It is enough to show that
$\g(R)\g(G(\mu)) \in \im(\k)$ for all $\mu \in \Lambda(n, r)$, using the
decomposition of the identity of $\sqnr$ into orthogonal idempotents
$\phi_{\mu, \mu}^1$.  Lemma 3.4.3 and Remark 3.4.4 show
that the element $\phi_{\l, \w}^1
\phi_{\w, \w}^w \phi_{\w, \mu}^1$ maps under $\k$ to a multiple of
$\g(R) \g(G(\mu))$, where $w = \rho^{-\mu_n}$ and where $\l$
is such that $R$ sends the weight space $\mu$ to the weight space
$\l$.  The case for $R^{-1}$ is similar but uses $w = \rho^{+\mu_1}$.
\qed\enddemo

The next result also makes use of specialization.  There is no problem
in replacing $v$ by $1$ here since we are working over an integral
form where $v-1$ cannot occur as a factor of a denominator.

\proclaim{Lemma 3.4.7}
The map $\k$, restricted to the ${\Bbb Z}[q, q^{-1}]$-form of $\sqnr$ given
in Theorem 2.2.4, is a monomorphism when the parameter $v$ is 
specialized to $1$.
\endproclaim

\demo{Proof}
Using an easy weight space argument, this can be reduced to showing
that the images of basis elements $\phi_{\l, \mu}^d$, for fixed $\l$
and $\mu$, are independent.

Suppose $$
k = \k\left( \sum_{d \in \D_{\l, \mu}} c_d \phi_{\l, \mu}^d \right) = 0
.$$  This means that $k$ acts as zero on $$
e_\mu := e_{\ell(\mu)_1} \otimes \cdots \otimes e_{\ell(\mu)_r}
.$$  

Consider the action of one $\phi_{\l, \mu}^d$ which appears in
the sum on $e_\mu$.  By using the factorisation of $\phi_{\l, \mu}^d$
as in the proof of Proposition 2.5.1 together with lemmas 3.4.3 and
3.4.5 (in the case $v = 1$), the isomorphism of Lemma 3.3.3 and the
definition of $\k$ in Lemma 3.4.2, we have $$
\phi_{\l, \mu}^d . e_\mu = 
\sum_{w \in d \W_\mu}
l_w e_{\l(w(1))} \otimes \cdots \otimes e_{\l(w(r))}
,$$ where $\l(i) := \ell(\l)_i$, and the scalars $l_w$ are positive
integers.  It follows from this that the actions of $\phi_{\l, \mu}^d$
and $\phi_{\l, \mu}^{d^\prime}$ on $e_{\mu}$
(where $d, d'$ lie in distinct double
$\W_\l$--$\W_\mu$-cosets) are nonzero and have distinct supports.  
The claim follows.
\qed\enddemo

We are now ready to show that the two definitions of the affine
$q$-Schur algebra are compatible with each other.

\proclaim{Theorem 3.4.8}
Over $\qv$, the map $\k$ gives an isomorphism between the 
affine $q$-Schur algebra $\sqnr$ and the quotient of
$\ugn$ by the kernel of its action on tensor space.
\endproclaim

\demo{Proof}
By Lemma 3.4.6, it is enough to show that $\k$ is a monomorphism.  To
do this, we consider an element in $\ker \k$ and then use the
techniques of Lemma 3.3.4 to show that unless this element is zero, we
have a counterexample in the case where $v$ is specialized to $1$,
which contradicts Lemma 3.4.7.
\qed\enddemo

\head 4. Concluding results \endhead

\subhead 4.1 Comparison of tensor space with $q$-tensor space \endsubhead

We now show that the $q$-tensor space of \S2.3 is isomorphic as an
$\sqnr$-module to tensor space $V^{\otimes r}$.  In \S4, we continue
to use $\qv$ as a base ring unless otherwise stated.

\proclaim{Lemma 4.1.1}
As a $\ugn$-module, tensor space $V^{\otimes r}$ is generated by any
vector $$
e_w = e_{w(1)} \otimes \cdots \otimes e_{w(r)}
$$ for $w \in \W$.
\endproclaim

\demo{Proof}
Using Lemma 3.4.3, it is enough to show that $V_\w$ is generated as a
$\ugn$-module by $e_w$.  We introduce endomorphisms
$\t(T_{s_i}^{-1})$ of $V_\w$ corresponding to the action of the
elements $v^{-1}E_i F_i - 1$ of $\ugn$ on $V_\w$.
It is easily checked that $\t(T_{s_i}^{-1})$ is the inverse of
$\t(T_{s_i})$.

Let $e_{j_1} \otimes \cdots \otimes e_{j_r}$ be a basis element of
$V_\w$.  Let $1 \leq p(j, i) \leq r$ be the unique integer such that
$j_{p(j, i)} \equiv i \mod n$.  It is not hard to show that if $p(j, i) <
p(j, i+1)$ then $$
\t(T_{s_i})(e_{j_1} \otimes \cdots \otimes e_{j_r})
=
v^* e_{s_i(j_1)} \otimes \cdots \otimes e_{s_i(j_r)}
,$$ for some integer $*$, 
and if $p(j, i) > p(j, i+1)$ (the only other possibility), we have $$
\t(T_{s_i}^{-1})(e_{j_1} \otimes \cdots \otimes e_{j_r})
=
v^* e_{s_i(j_1)} \otimes \cdots \otimes e_{s_i(j_r)}
.$$  If $w = \rho^{\pm 1}$, we have $$
\t(T_w)(e_{j_1} \otimes \cdots \otimes e_{j_r})
=
v^* e_{w(j_1)} \otimes \cdots \otimes e_{w(j_r)}
.$$

We can now use the techniques of Lemma 3.3.3 to deduce the claim.
\qed\enddemo

\proclaim{Theorem 4.1.2}
As $\sqnr$-modules, $q$-tensor space and tensor space are isomorphic.
The isomorphism may be chosen to take $\phi_{\w, \w}^1$ to $$
e_\w = e_1 \otimes \cdots \otimes e_r
.$$
\endproclaim

\demo{Proof}
By Lemma 4.1.1, we need only show that no element of $q$-tensor space
annihilates the vector $e_\w$.  Suppose $$
k = \sum_{\l \in \Lambda(n, r), d \in \D_\l} c_{\l, d} \phi_{\l, \w}^d
$$ annihilates $e_\w$.  We may assume that only one $\l$ occurs in the
sum.  In this case, $$
k = \phi_{\l, \w}^1 . 
\left( \sum_{d \in {\D_\l}} c_{\l, d} \phi_{\w, \w}^d \right)
,$$ and if $k$ acts as zero, so does $$
\sum_{w \in \W_\l} \sum_{d \in {\D_\l}} c_{\l, d} \phi_{\w, \w}^{wd}
.$$  (To see the latter assertion, premultiply by $\phi_{\w, \l}^1$
and use Lemma 3.4.5.)  This contradicts Lemma 3.3.4 unless all the
$c_{\l, d}$ are zero.
\qed\enddemo

\remark{Remark 4.1.3}
The isomorphism between $q$-tensor space and tensor space is not a
natural one as in the case of the ordinary $q$-Schur algebra or the
case $q = 1$.  To see this, compare the actions of $\phi_{\w,
\w}^{s_r}$ on the generating vector: in the case of $q$-tensor space,
this produces a combination of one other basis vector, but in the case
of tensor space, $e_1$ occurs to the left of $e_r$ in $e_\w$,
so a combination of two basis vectors is produced.
\endremark

\subhead 4.2 Presentations for the affine Hecke algebra \endsubhead

An important consequence of Theorem 4.1.2 is that it provides, in
conjunction with Theorem 2.3.3, a right action of $\H$ on tensor
space, and furthermore, the centralizing algebra of the action of
$\ugn$ on tensor space is isomorphic to $\H$.  This means that the
endomorphisms $y_i$ of Definition 3.2.4 can be regarded as elements of $\H$.

The results of this section consider this action further to show
that our affine Hecke algebra is in fact
isomorphic to the affine Hecke algebra in the sense of \cite{{\bf 2},
3.1 Definition}.

\proclaim{Lemma 4.2.1}
Let $$
e_{j_1} \otimes \cdots \otimes e_{j_r}
$$ be a basis element of $V_\w$ satisfying $1 \leq j_t \leq r$ for all
$t$.  Suppose that $j_i = r-1$ and $j_{i+1} = r$.  Then we have $$
(e_{j_1} \otimes \cdots \otimes e_{j_r}) . T_{s_i} y_i T_{s_i} = 
(e_{j_1} \otimes \cdots \otimes e_{j_r}) . v^2 y_{i+1}
.
$$
\endproclaim

\demo{Proof}
By definition of the endomorphism $y_{i+1}$ in Definition 3.2.4, the
right hand side is easily seen to be $$
v^2 e_{j_1} \otimes \cdots \otimes e_{j_i} \otimes e_{-n + j_{i+1}} 
\otimes e_{j_{i+2}} \otimes \cdots \otimes e_{j_r}
.$$  Because of the way $S_q(n, r)$ embeds in $\sqnr$, the action of
$T_{s_i}$ on $V_n^{\otimes r}$ is the same as in the finite case.  We
therefore use the definition of this action in \cite{{\bf 9}, \S1.5} to
show that  $$
(e_{j_1} \otimes \cdots \otimes e_{j_r}) . T_{s_i} y_i = 
v e_{j_1} \otimes \cdots \otimes 
e_{-n + j_{i+1}} \otimes e_{j_i} 
\otimes \cdots \otimes e_{j_r}.
$$  It suffices to show that $$
(e_{j_1} \otimes \cdots \otimes e_{-n + j_{i+1}} \otimes e_{j_i} 
\otimes \cdots \otimes e_{j_r}) . T_{s_i} = 
v (e_{j_1} \otimes \cdots \otimes e_{j_i} \otimes e_{-n + j_{i+1}} 
\otimes \cdots \otimes e_{j_r}).
$$  Recall that $-n + j_{i+1} = -n + r$.  There are no other indices
with values between $-n + r - 1$ and $0$.  

It is enough to establish $$
{(e_{(n-r+1) + j_1} \otimes \cdots \otimes e_{1 - r + j_{i+1}} \otimes
e_{(n-r + 1) + j_i} \otimes \cdots \otimes e_{(n-r+1) + j_r})
. T_{s_i} \quad \quad
\hfill}\atop
{\hfill = v (e_{(n-r+1) + j_1} \otimes \cdots \otimes e_{(n-r+1) + j_i} \otimes
e_{1 - r + j_{i+1}} \otimes \cdots \otimes e_{(n-r+1) + j_r}),}
$$ and then act $$
R^{- n + r - 1} \in \ugn
$$ on the left of both sides.
Since all the indices lie between $1$ and $n$ and $$n = (n-r+1) + j_i >
1 - r + j_{i+1} = 1,$$ the claim follows from the action of $\H$ on 
$V_n^{\otimes r}$.
\qed\enddemo

\proclaim{Corollary 4.2.2}
We have $T_{s_i} y_i T_{s_i} = v^2 y_{i+1}$.
\endproclaim

\demo{Proof}
This follows from Lemma 4.1.1, Lemma 4.2.1 and the double centralizer
property for $\ugn$ and $\H$ on tensor space.
\qed\enddemo

\proclaim{Lemma 4.2.3}
Suppose $1 \leq i < r$, $1 \leq j \leq r$ and $j \ne i, i+1$.
Let $$
e_{k_1} \otimes \cdots \otimes e_{k_r}
$$ be a basis element of $V_\w$ satisfying $1 \leq k_t \leq r$ for all
$t$.  Suppose that $k_i = r-2$, $k_{i+1} = r-1$ and $k_j = r$.
Then we have $$
(e_{k_1} \otimes \cdots \otimes e_{k_r}) . y_j T_{s_i} = 
(e_{k_1} \otimes \cdots \otimes e_{k_r}) . T_{s_i} y_j.
$$
\endproclaim

\demo{Proof}
This is proved using an argument very similar to the that in 
proof of Lemma 4.2.1.
\qed\enddemo

\proclaim{Corollary 4.2.4}
We have $y_j T_{s_i} = T_{s_i} y_j$ if $j \ne i, i+1$.
\endproclaim

\demo{Proof}
This uses Lemma 4.1.1, Lemma 4.2.3 and the double centralizer property.
\qed\enddemo

We can now establish an isomorphism between $\H$ and the affine Hecke
algebra of \cite{{\bf 2}, \S3}.

\proclaim{Theorem 4.2.5}
The algebra $\H(\W)$ is isomorphic to the algebra given by generators $$
\{\s_i^{\pm 1} : 1 \leq i < r\} \cup
\{y_i^{\pm 1} : 1 \leq i \leq r\}
$$ and defining relations $$\eqalignno{
\s_i \s_i^{-1} &= \s_i^{-1} \s_i = 1, & (1)\cr
\s_i \s_{i+1} \s_i &= \s_{i+1} \s_i \s_{i+1}, & (2)\cr
\s_i \s_j &= \s_j \s_i \text{\ if } |i - j| > 1, & (3)\cr
(\s_i + 1)(\s_i - v^2) &= 0, & (4)\cr
y_j y_j^{-1} &= y_j^{-1} y_j = 1, & (5)\cr
y_j y_k &= y_k y_j, & (6)\cr
y_j \s_i &= \s_i y_j \text{\ if } j \ne i, i+1, & (7)\cr
\s_i y_i \s_i &= v^2 y_{i+1}. & (8)\cr
}$$ The isomorphism may be chosen to identify $\s_i$ with $T_{s_i}$.
\endproclaim

\demo{Proof}
We identify $\s_i$ with $T_{s_i}$.  Relations (1)--(4) are known from
Definition 1.2.1.  Relations (5) and (6) are obvious from considering
the faithful action of $\H$ on tensor space.  Relations (7) and (8)
follow from corollaries 4.2.4 and 4.2.2 respectively.

These relations are sufficient since \cite{{\bf 2}, 3.1 Lemma} shows that a
basis for the algebra given by these generators and relations has a
basis $$
y_1^{c_1} \cdots y_r^{c_r} T_w
$$ as $c_i \in {\Bbb Z}$ and $T_w$ ranges over a basis for $\H({\Cal
S}_r)$.  

Since the algebra in the statement is clearly isomorphic to its own
opposite algebra, there is also a basis of the form $$
T_w y_1^{c_1} \cdots y_r^{c_r}.
$$  It is not hard to verify that the map $$
T_w y_1^{c_1} \cdots y_r^{c_r} \mapsto 
e_\w . T_w y_1^{c_1} \cdots y_r^{c_r}
$$ is injective, from which we deduce that the algebra
of \cite{{\bf 2}} acts faithfully on tensor space under the identifications
given.  It follows that the relations given in the statement are sufficient.
\qed\enddemo

\remark{Remark 4.2.6}
Using Theorem 4.2.5 we can give a third construction of the affine
$q$-Schur algebra more compatible with the results of \cite{{\bf 2}}.  This
is formally the same as Definition 2.1.4, but uses the version of the
affine Hecke algebra given in \cite{{\bf 2}}.
\endremark

The presentation for $\H$ just given is more convenient when dealing
with tensor space, in contrast to Proposition 1.2.3 which is adapted
for $q$-tensor space.  Using the basis $T_w \Pi y_i^{c_i}$ of $\H$ which
was mentioned in the proof of the theorem, one can explicitly
calculate the action of the generators $T_{s_i}$ on tensor space.  We
do not do this here, because it is much more complicated than in the
case of type finite $A$.  In type finite $A$, the action of $T_{s_i}$
on a basis element of $V_n^{\otimes r}$ produces a linear combination of at
most two basis elements, but in the case of the 
action of $T_{s_i}$ on $V^{\otimes r}$ in type affine $A$, it is possible for
arbitrarily many basis elements to be involved.  This happens because
we are using a basis of $\H$ whose structure constants are much more
complicated than those of the basis given in Proposition 1.2.3.  (This
complication arises from relation (8) in the statement of the theorem.)

\subhead 4.3 Some questions \endsubhead

It seems from the results of \S3 that the Hopf algebra $\ugn$ is a
natural object.  Further evidence for this is provided by the results
of \cite{{\bf 25}, \S11}, where it is proved that the affine $q$-Schur
algebra introduced in this paper occurs naturally in a $K$-theoretic
context.  This generalizes the interpretation in \cite{{\bf 1}} of 
the ordinary $q$-Schur algebra
in terms of the geometry of the action of the general linear group
(over a finite field) acting on certain pairs of flags.

Comparing the situation with the relationship between
$\ugln$ and the ordinary $q$-Schur algebra raises some more questions
about $\ugn$, such as the following ones.

\proclaim{Question 4.3.1}
Can one describe a suitable $\A$-form of the algebra $\ugn$ which maps under
$\a_r$ to the natural $\A$-form of $\sqnr$?
\endproclaim

This is an important question because it would show that our results
also hold at roots of unity.

\proclaim{Question 4.3.2}
To what extent do these results hold in the case where $n < r$?
\endproclaim

In this case, there is no weight $\w$, and
one should not expect the Hecke algebra to act
faithfully on tensor space, since it does not act faithfully in the case of the
ordinary $q$-Schur algebra.  In the latter case, it can be shown that
the centralizing algebra of the $q$-Schur algebra is a quotient of the
Hecke algebra.  Something similar may hold in the affine case.

A main obstacle to proving this lies in the difficulty of embedding
$\sqnr$ in $\widehat S_q(n+1, r)$.  This is
in contrast to the finite case, where there are obvious embeddings.
The lack of an obvious determinant map also means that it is not clear
how to realize $\sqnr$ as a quotient of $\widehat S_q(n, r+n)$.

\proclaim{Question 4.3.3}
Is there a basis for $\ugn$ similar to Lusztig's canonical basis for
the modified quantized enveloping algebras $\dot U$ which is compatible
with the Kazhdan--Lusztig basis of \S2.4 via the maps $\a_r$?
\endproclaim

Lusztig's algebras $\dot U$ are defined in \cite{{\bf 23}, \S23}.  They are
equipped with canonical bases $\dot B$, generalizing the canonical
bases for the algebras $U^+$ (which are the subalgebras of quantized
enveloping algebras generated by the elements $E_i$).  In the case of
the ordinary $q$-Schur
algebra, these bases have been shown in \cite{{\bf 10}} to be compatible
with the Kazhdan--Lusztig bases for the $q$-Schur algebra, in the
sense that the basis elements $\dot B$ map either to zero or to
Kazhdan--Lusztig basis elements.  It seems reasonable to hope
therefore that an analogous result exists for the affine case.

\proclaim{Question 4.3.4}
Do the structure constants for the affine $q$-Schur algebra with
respect to the Kazhdan--Lusztig basis lie in ${\Bbb N}[v, v^{-1}]$?
\endproclaim

This property is true for the ordinary $q$-Schur algebra, as was
proved in \cite{{\bf 14}}.  If the more general result is true, it
should have an interpretation in terms of intersection cohomology.

\head Acknowledgements \endhead

The author thanks J. Du and R.J. Marsh for helpful conversations.
\vfill\eject

 \vfill\eject
\end